\documentclass[reqno]{amsart}

\input xy
\xyoption{all}
\usepackage{accents}
\usepackage{epsfig}
\usepackage{color}
\usepackage{amsthm}
\usepackage{amssymb}
\usepackage{amsmath}
\usepackage{amscd}
\usepackage{amsopn}
\usepackage{graphicx}

\usepackage{xspace}

\usepackage{url}
\usepackage{enumitem, hyperref}\hypersetup{colorlinks}

\usepackage{color} %\textcolor{red}{Test}

\definecolor{darkred}{rgb}{1,0,0} %can change the intensity in [0,1]
\definecolor{darkgreen}{rgb}{0,0.8,0}
\definecolor{darkblue}{rgb}{0,0,1}

\hypersetup{colorlinks,
linkcolor=darkblue,
filecolor=darkgreen,
urlcolor=darkred,
citecolor=darkgreen}

\makeatletter
\def\reflb#1#2{\begingroup
    #2%
    \def\@currentlabel{#2}%
    \phantomsection\label{#1}\endgroup
}
\makeatother

\numberwithin{equation}{section}
\numberwithin{Theorem}{section}

\theoremstyle{definition}

\theoremstyle{remark}

\def    \eps    {\epsilon}

\newcommand{\CA}{{\mathcal A}}

\newcommand{\CS}{{\mathcal S}}

\newcommand{\fc}{{\mathfrak c}}

\def    \R      {{\mathbb R}}

\def    \Z      {{\mathbb Z}}

\def    \Q      {{\mathbb Q}}

\def    \12    {{\frac{1}{2}}}

\def    \p      {\partial}
\def    \sel  {\operatorname{c}}

\def    \HF     {\operatorname{HF}}
\def    \HC     {\operatorname{HC}}
\def    \CC     {\operatorname{CC}}

\def    \H     {\operatorname{H}}

\def    \MUCZ  {\operatorname{\mu_{\scriptscriptstyle{CZ}}}}

\begin{document}

%%%%%%%%%%%%%%%%%%%%%%%%%%%%%%
%   TEXT FORMATTING

\setlength{\smallskipamount}{6pt}
\setlength{\medskipamount}{10pt}
\setlength{\bigskipamount}{16pt}

%%%%%%%%%%%%%%%%%%%%%%%%%%

%%%%%%%%%%%%%%%%%%%%%%%%%%

%%%%%%%%%%%           BEGINNING OF  TEXT

%%%%%%%%%%%%%%%%%%%%%%%%%%

\title[Shopping List]{My Contact Homology Shopping List}

\author[Viktor Ginzburg]{Viktor L. Ginzburg}
%\author[Ba\c sak G\"urel]{Ba\c sak Z. G\"urel}

% \address{BG: Department of Mathematics, University of Central Florida,
%   Orlando, FL 32816, USA} \email{basak.gurel@ucf.edu}

\address{Department of Mathematics, UC Santa Cruz, Santa Cruz, CA
  95064, USA} \email{ginzburg@ucsc.edu}



\date{\today} 

\thanks{The work is partially supported by NSF grant DMS-1308501}

\bigskip

\begin{abstract}
  We list the properties of contact homology, beyond purely formal,
  needed for the proofs of some of the recent applications of contact
  homology in dynamics to work. The list is put together for the AIM
  Transversality in Contact Homology Workshop.

\end{abstract}

\maketitle

\tableofcontents

\section{Introduction}

In this write-up for the AIM Transversality in Contact Homology
Workshop, we identify the properties of contact homology, beyond
purely formal, needed for the proofs in \cite{GGM,GHHM} to work. We
also expect these properties to be useful for other applications of
contact homology in dynamics; see, e.g., \cite{Gu:pr,HM}.

The main result of \cite{GHHM} is the ``SDM theorem'' asserting that
under natural additional conditions the presence of a simple closed
Reeb orbit of a particular type (the so-called SDM) implies the
existence of infinitely many simple closed Reeb orbits. As an
application, we (re)prove the existence of at least two closed Reeb
orbits for any contact form supporting the standard contact structure
on $S^3$; see also \cite{CGH,GGo,LL} for other proofs of this or for more
general results. In \cite{GGM}, we prove a variant of the Conley
conjecture for Reeb flows on the pre-quantization circle bundles over
aspherical symplectic manifolds and then use this fact to establish
the existence of infinitely many simple closed orbits for a low energy
twisted geodesic flow on a surface of positive genus with
non-vanishing magnetic field. We refer the reader to \cite[Sect.\ 4
and 5]{GG:CCsurvey} for a much more detailed discussion of these
results.

In \cite{GHHM}, we work exclusive with the linearized contact
homology, and hence one should be able to reprove the results of that
paper relying instead on the equivariant symplectic homology and thus
bypassing the transversality problems; cf.\ \cite{BO:eq}. (Such as
translation of the proof to a different language is however rather
non-trivial.) On the other hand, in \cite{GGM}, it is essential to use
mainly the cylindrical contact homology with its additional grading by
the free homotopy classes of loops. (The SDM theorem from \cite{GHHM}
also enters the proof as an ingredient. It is interesting to note
that, as of this writing, the SDM theorem does not have, due to the
index restrictions, an analog relying on the cylindrical contact
homology. Working with twisted geodesic flows, one can circumvent the
use of cylindrical contact homology by utilizing a very particular
filling.)

Throughout this write-up, we will focus on the cylindrical contact
homology. There are two reasons for this. First of all, this is the
type of contact homology where it is not clear how to get around the
transversality problems. Secondly, much of what we say readily translates to
the linearized contact homology.

It might be worth pointing out that for the (finite energy)
holomorphic curves, i.e., when the transversality problem is left
aside, the properties we list below are either well known or seem to
present no serious difficulty to prove. However, it is not \emph{a priori}
clear that these properties would automatically carry over once the
transversality problem is resolved and holomorphic curves are possibly
replaced by some other objects. Finally, we note that our approach to
the construction of the local contact homology is somewhat different,
at least on the technical level, from that in~\cite{HM}.

\medskip
\noindent\textbf{Acknowledgements.} 
The author is grateful to Daniel Cristofaro-Gardiner, Joel Fish,
Ba\c{s}ak G\"urel, Umberto Hryniewicz, Michael Hutchings, Eleny Ionel,
Leonardo Macarini, Gordana Matic, Otto van Koert and Katrin Wehrheim
for useful remarks and discussions.

\section{Contact homology}
\label{sec:ch}

\subsection{Setting}
In what follows, $(M^{2n-1},\xi)$ is a closed contact manifold. Let
$\alpha$ be a non-degenerate contact form with $\ker \alpha=\xi$. We
assume that the periodic orbits of $\alpha$ meet the necessary index
conditions for the linearized contact homology $\HC_*(M,\xi)$ and the
action-filtered linearized contact homology $\HC_*^{I}(M,\alpha)$,
where $I=(a,\,b)$, to be defined. (We will always require $a$ and $b$
to be outside the action spectrum $\CS(\alpha)$ of $\alpha$.)  For the
sake of simplicity, let us also assume that $c_1(\xi)=0$ in $\H^2(M;\Z)$.

The (cylindrical) contact homology complex $\CC_*(M,\alpha)$ is
generated by the good orbits $x$ of $\alpha$, graded by
$|x|=\MUCZ(x)+n-3$ and filtered by the action. The differential
$\p\colon \CC_*(M,\alpha)\to \CC_{*-1}(M,\alpha)$ depends on some
auxiliary data and has the form
$$
\p x=\sum_y m(x,y)\,y,
$$
where $m(x,y)$ counts certain maps $u\colon S^1\times \R\to M$ or some
other geometric objects with signs and weights. In what follows, for
the sake of simplicity, we assume that $u$ is the $M$-component of a
map $S^1\times \R\to M\times \R$ which must satisfy some (translation
invariant) equation or more generally certain conditions depending on
some auxiliary structures and choices, which we
actually do not quite know at this moment and which we call (CR)
here. In the original definition of the contact homology, (CR) is the
Cauchy--Riemann equation in the symplectization and, in addition, the
solutions are required to have finite Hofer energy; see \cite{SFT} and
also \cite{Bo}. For the sake of brevity, let us refer, somewhat
unconventionally, to $u$ as a Floer trajectory. Note that (CR) should
make sense even when $\alpha$ is degenerate although  one can probably
circumvent this requirement.

The $\omega$-energy of $u$ is by definition
$$
E_\omega(u)=\int_u d\alpha.
$$
Clearly, the $\omega$-energy is translation invariant, and one
essential feature of (CR) should be that
\begin{equation}
\label{eq:zero-energy}
E_\omega(u)=0\textrm{ iff  $u$ is a trivial Floer
trajectory},
\end{equation}
i.e., a closed Reeb orbit. (See, e.g., \cite[Lemma 5.4]{SFT:comp} for
a proof of \eqref{eq:zero-energy} in the holomorphic case.)
Furthermore,
$$
E(u)=\CA_\alpha(x)-\CA_\alpha(y),
$$
where $u$ is (partially) asymptotic to $x$ and $y$.  Here the action
$\CA_\alpha$ is defined as
$$
\CA_\alpha(x):=\int_x\alpha.
$$
Thus $\p$ is action decreasing. The complex $\CC_*(M,\alpha)$, and
hence the homology, is also graded by the free homotopy classes of
loops in $M$. The filtered contact homology is defined, by continuity,
even when $\alpha$ is degenerate, provided that $\alpha$ has
perturbations meeting the index conditions; see Section \ref{sec:gen}.

\subsection{Spatial localization vs.\ energy localization} 
\label{sec:mon}
In this section, we do not require $\alpha$ to be non-degenerate or
$M$ to be compact. There may be a free homotopy class of loops, say
$\fc$, or a collection of such classes fixed in the background.

In several instances we need to have lower bounds on the
$\omega$-energy of Floer trajectories passing through a certain region
to spatially localize low energy trajectories, i.e., to ensure that
such trajectories are confined to the complement of the region. Here
is one variant of such an assertion, which would probably cover all
the instances where the localization has been used so far.

Let $S$ be a compact subset of $M$ and let $A$ be the set of actions
of closed Reeb orbits passing through $S$. (The set $S$ is usually a
hypersurface in $M$.)  Fix a neighborhood $N$ of $S$.  Then we want to
be able to say that there exists a constant $\eps=\eps(S,N,\alpha)>0$
such that
\begin{equation}
\label{eq:elb}
E_\omega(u)>\eps
\end{equation}
for every $u$ passing through $S$, provided that $u$ is (partially)
asymptotic to some $x$ and $y$ which are contained entirely in the
complement of $N$ and such that the action interval
$[\CA_\alpha(y),\CA_\alpha(x)]$ does not intersect $A$.

Moreover, $\eps>0$ can be taken so that %this is true
\eqref{eq:elb} holds for every Floer trajectory for every contact form
$\alpha'$ which is $C^\infty$-close to $\alpha$. (We probably don't
care if we exactly have $\ker\alpha'=\ker\alpha$ or not.) It would
also be useful to know that $\eps$ depends only on $\alpha|_N$. This
is, of course, a rather general assertion, and what is actually needed
is its various particular cases as in, e.g., \cite{GHHM}.  When Floer
trajectories come from genuine holomorphic curves in the
symplectization, \eqref{eq:elb} follows from, e.g., the variant of
Gromov compactness proved in \cite{Fi}, similarly to an argument in
\cite{McL}, combined with \eqref{eq:zero-energy}.

Here is a typical application.  Assume that $x$ and $y$, closed Reeb
orbits of $\alpha'$, lie in a complement of $N$ and that the action
difference $|\CA_{\alpha'}(x)-\CA_{\alpha'}(y)|$ is small. Then there
are no Floer trajectories between $x$ and $y$ passing through $S$
provided that there are no closed Reeb orbits of $\alpha$ in the free
homotopy class of $x$ and $y$ passing through $S$, i.e.,
$A=\emptyset$. This argument can be used to show that $\p^2=0$ in the
construction of local contact homology.  Note, however, that the
latter condition appears to be very difficult to check in general
(except perhaps in very few cases such as completely integrable Reeb
flows) unless there are non-trivial homotopy class or action
restrictions. A variant of \eqref{eq:elb} is utilized is the proof of the
SDM theorem in \cite{GHHM}.

\subsection{The shift map and the Lusternik--Schnirelmann theory}
The cylindrical contact homology has a degree-two downward (i.e.,
degree $-2$) shift map $D$ which is an analogue of the pairing with
the generator of $\H^2(BS^1;\Z)$ in the equivariant symplectic
homology; see \cite[Sect. 7.2]{BO}. It would be useful to know how $D$
effects the spectral invariants of $\alpha$. More specifically, denote
by $\sel_w(\alpha)$ the spectral invariant of $\alpha$ associated with
$w\in\HC_*(M,\xi)$. (By definition, we set $\sel_0=-\infty$.) Note
that although we are not making any non-degeneracy assumptions,
$\sel_w(\alpha)$ is defined on the level of homology. It should be
very much straightforward to show that, for any reasonable definition
of $D$, we automatically have $\sel_{D(w)}(\alpha)\leq
\sel_w(\alpha)$.  However, what we want is the strict inequality,
i.e., as in other versions of the Lusternik--Schnirelmann (LS) theory
(cf.\ \cite[Sect.\ 6]{GG:gaps}),
\begin{equation}
\label{eq:LS}
\sel_{D(w)}(\alpha)<\sel_w(\alpha)
\end{equation}
when $w\neq 0$ and all simple closed Reeb orbits of $\alpha$ are
isolated in the extended phase space. The latter condition is
essential. In the non-degenerate case, \eqref{eq:LS} is a consequence
of \eqref{eq:zero-energy} once we know that $D$ can be defined via
counting solutions of the (CR) equation. Without non-degeneracy, one
should in addition show that, roughly speaking, $D=0$ on the level of
local contact homology. A variant of \eqref{eq:LS} for the ECH is used
in, e.g., \cite{CGH}; see also \cite{Gu:pr} for an application of the
shift map in the context of equivariant symplectic or linearized
contact homology. I am not aware of any situation where having
\eqref{eq:LS} for cylindrical contact homology is crucial, but it is
probably only a matter of time before we encounter one.

\section{Cobordisms and continuation}
\subsection{Generalities}
\label{sec:gen}
In this section, a contact manifold is a manifold equipped with a
contact form rather than a contact structure.  Let us assume that, as
expected, a symplectic cobordism $(V,\omega)$ from a closed contact
manifold $(M_1,\alpha_1)$ to another one $(M_0,\alpha_0)$ induces a
``continuation map'' between filtered contact homology complexes when
the forms are non-degenerate. The map on the level of the filtered
homology is always defined when the end points of the action are
outside $\CS(\alpha_1)$ and $\CS(\alpha_0)$, regardless of whether the
forms are non-degenerate or not. The continuation maps are again
obtained by counting (finite energy) solutions of a certain equation
in $\hat{V}$ obtained from $V$ by
attaching the bottom and the top parts of the symplectizations of
$M_1$ and, respectively, $M_0$. This equation, which we
call (CR) again, depends on some extra structure, etc. Note that at this
point we do not quite know what (CR) should be.

Doing dynamics, it is sometimes more convenient to think, by analogy
with the Hamiltonian setting, in terms of monotone homotopies. Here is
a formal definition.  We say that a family of contact forms
$\alpha_s$, $s\in [0,\,1]$, on $M$ foliates a symplectic cobordism
$(V,\omega)$ from $\alpha_1$ to $\alpha_0$ if there exists a family of
contact type embeddings $j_s\colon M\to V$ smoothly foliating $V$, in
the obvious sense, and such that $j_s^*\omega=d\alpha_s$ on $M$. Note
that then $V$ is necessarily diffeomorphic to the cylinder $M\times
[0,\,1]$ with $\p V$ comprising two parts: $j_0(M)$ (the positive end)
and $j_1(M)$ (the negative end).  Thus, with this convention, the
family $\alpha_s$ is ``decreasing'' and we have a map
$\HC_*^I(\alpha_0)\to \HC_*^I(\alpha_1)$ when the end points of $I$
are outside $\CS(\alpha_0)\cup\CS(\alpha_1)$.

Then we have the following ``stability result'' for contact homology:
the map $\HC_*^I(\alpha_0)\to \HC_*^I(\alpha_1)$ is an isomorphism
when the end points of $I$ are outside $\CS(\alpha_s)$ for all $s$.
Furthermore, let $I_s$ be a family of intervals such that for every
$s$ the end points of $I_s$ are outside $S(\alpha_s)$ for a family
$\alpha_s$ not necessarily foliating a cobordism. Then the contact
homology spaces $\HC_*^{I_s}(\alpha_s)$ are isomorphic. In particular,
the global contact homology depends only on the contact structure, the
filtered contact homology is defined for degenerate forms, etc. This
result is stated (without a proof and for the linearized contact
homology) in \cite{GHHM}. The proof differs from the Hamiltonian case,
where a homotopy induces a map in both directions, but it is
nonetheless almost entirely formal, up to one property of the
continuation maps, and hence carries over to other types of contact
homology theories.

This property is that when $\alpha$ is non-degenerate and $|\tau|$
is sufficiently small there should be a matching choice of the auxiliary
data for the forms $\alpha$ and $(1+\tau)\alpha$ resulting in an
isomorphism between their contact homology complexes. On the level of
the filtered contact homology this isomorphism should be equal to the
one induced by the natural cobordism between these two forms.

When dealing with the cylindrical homology, one has to assume in
addition that, say, for a dense set of $s\in [0,\,1]$ the forms
$\alpha_s$ foliating a cobordism have non-degenerate perturbations
meeting the index conditions.

\subsection{Spatial localization for continuation maps} In some cases,
it is essential to know that low energy solutions of (CR) in $\hat{V}$
cannot pass through a certain region in $M$. To be more specific,
assume that a family $\alpha_s$ foliates $V$. Furthermore, let, as in
Section \ref{sec:mon}, $S$ be a compact subset of $M$ and $N$ be a
neighborhood of $S$. Denote by $A$ the union of the action spectra of
$\alpha_0$ and $\alpha_1$ for closed Reeb orbits passing through $S$.
Then we need to know that there exists a constant
$\eps=\eps(S,N,\alpha_s)>0$ such that for every continuation
trajectory $u$ passing through $S$ and (partially) asymptotic to some
$x$ and $y$ we have
\begin{equation}
\label{eq:elb2}
\CA_{\alpha_0}(x)-\CA_{\alpha_1}(y)
>\eps,
\end{equation}
provided that $x$ and $y$ are contained entirely in the complement of
$N$ and that the action interval
$[\CA_{\alpha_1}(y),\CA_{\alpha_0}(x)]$ does not intersect $A$.

One point to keep in mind here is that the the auxiliary structure on
$\hat{V}$ we use here, and hence (CR), are adapted to the family
$\alpha_s$.

Above, we may again have a free homotopy class of loops in $M$ or a
collection of such classes fixed in the background. Moreover, we would
need to know that $\eps>0$ can be taken so that \eqref{eq:elb2} still
holds when the family $\alpha_s$ is replaced by another foliating
family of contact forms $\alpha'_s$ which is $C^\infty$-close to
$\alpha_s$, and, ideally, that $\eps$ depends only on $\alpha_s|_N$.

Of course, what we really need are some particular cases of
\eqref{eq:elb2}. In \cite{GHHM}, we used the case of
\eqref{eq:elb2} (for the linearized contact homology) with
$N=S^1\times N_0$, where $N_0$ is a spherical shell in $\R^{2n}$ and
$\alpha_s=\lambda+c(s)\,dt$ and $S=S^1\times S^{2n-1}$. Here $\lambda$
is a primitive of the standard symplectic structure on $\R^{2n}$,
$c(s)$ is a decreasing family of constants, and $t$ is a coordinate on
$S^1$. Then \eqref{eq:elb2} holds for holomorphic curves when the
almost complex structure on $V$ comes from the standard complex
structure on $N_0$; see \cite{GHHM}.

Another instance where a variant of \eqref{eq:elb2} can be used is the
invariance of the local contact homology. In this
setting, $N$ is a thickened boundary of an isolating tubular
neighborhood of a closed Reeb orbit $x$. We would prefer to set
$\alpha_s=\alpha$ on $N$, but a constant family does not foliate a
cobordism. Instead, one can take a small perturbation of the constant
family: $\alpha_s=c(s)\alpha$, where $c(s)$ is a $C^2$-small monotone
decreasing function on $[0,\,1]$. (The argument in \cite{HM} is different.)

It is worth pointing out that \eqref{eq:elb2} has probably never been
proved in detail in the general form as stated above for holomorphic
curves, but from the first glance it looks correct.

\section{Local contact homology}
\subsection{Definitions and basic properties}
The local contact homology is associated to an isolated (in the
extended phase space) closed Reeb orbit. To be more specific, let $x$
be such an orbit of the Reeb flow of $\alpha$. We do not assume that
$x$ is non-degenerate or simple. Under a small non-degenerate
perturbation of $\alpha$ the orbit $x$ splits into a finite collection
of non-degenerate orbits $x'_i$ contained in a small isolating
neighborhood $U$ of $x$, and one can form a contact homology complex
using $x'_i$ with the differential defined exactly as in the global
case. We denote the resulting homology by $\HC_*(x)$ or
$\HC_*(\alpha,x)$. The local contact homology is invariant under
deformations $\alpha_s$ of $\alpha$ as long as $x$ is uniformly
isolated, i.e., having a common isolating neighborhood for all forms
$\alpha_s$. The local contact homology is discussed in detail in
\cite{HM} and then, in lesser detail, in \cite{GHHM}. Our outline
of the proof that this homology is defined and well-defined differs in
a number of ways from the proofs in \cite{HM}.

When the orbit $x$ is simple, there are no transversality problems in
the definition of the local contact homology because all orbits in
question belong to a simple free homotopy class. In this case, we have
\begin{equation}
\label{eq:lch-lfh1}
\HC_{*+n-3}(x)=\HF_*(\varphi),
\end{equation}
where on the right we have the local Floer homology of (the germ of)
the Poincar\'e return map $\varphi$ of $x$; see, e.g., \cite{GG:gaps}
for the definition. A variant of this identity is essentially
contained already in \cite[Sect.\ 6] {EKP}; see also \cite{HM} for
another proof. Note that in this case one can work over any
coefficient ring.

When $x$ is iterated, \eqref{eq:lch-lfh1} is no longer true. In this
case, the right hand side should be replaced by the equivariant
symplectic homology. To be more precise, assume $x=z^k$ where $z$ is
simple. Then the Poincar\'e return map $\varphi$ of $x$ is the $k$th
iteration of the Poincar\'e return map of $z$ and the
$\Z_k$-equivariant local Floer homology $\HF_*^{\Z_k}(\varphi)$ of
$\varphi$ is defined. We should have
\begin{equation}
\label{eq:lch-lfh2}
\HC_{*+n-3}(x)=\HF_*^{\Z_k}(\varphi),
\end{equation}
where both groups are taken over $\Q$; cf.\ \cite{BO:eq}. Now the
definition of the left hand side encounters the usual transversality
problems, and the definition of $\HC_*(x)$ and the proof of
\eqref{eq:lch-lfh2} should rely on a version of abstract
perturbations.

Let us ignore for a moment the transversality issue. Thus we assume
that there exists an almost complex structure on the symplectization
of $U$ meeting all the regularity assumptions. This is, roughly
speaking, equivalent to assuming that there exists a one-periodic in time
almost complex structure on a local cross section to $x$, which is regular for the
$k$-periodic map $\varphi$. In this case, as is observed in \cite{GHHM},
one can replace the equivariant homology $\HF_*^{\Z_k}(\varphi)$ by
the $\Z_k$-invariant homology $\HF_*(\varphi)^{\Z_k}$, and the equality 
\begin{equation}
\label{eq:lch-lfh3}
\HC_{*+n-3}(x)=\HF_*(\varphi)^{\Z_k}
\end{equation}
can then be proved along with \eqref{eq:lch-lfh2} by
reasonably conventional methods as in, e.g., \cite{EKP}. (The argument
is outlined in \cite{GHHM}.)

Overall, the situation here is rather similar to equating, as in
\cite{BO:eq}, the linearized contact homology and the equivariant
symplectic homology, and perhaps is even a bit simpler than that. Some
parts of the proof are given in \cite{GHHM}.

Although the identification \eqref{eq:lch-lfh2} is illuminating, it
has not been used anywhere to the best of the author's knowledge. What
has been used is a weaker result that $\dim \HC_*(x^k)$ is a bounded
function of $k$ as long as $x^k$ remains isolated; see \cite{HM}.

Finally, the local contact homology are the building blocks for the
global and filtered contact homology. Namely, assume that the only
point of the action spectrum $\CA(\alpha)$ in the interval $I$ is $c$
and that all orbits $x$ with action equal to $c$ are isolated. Then we
should have
\begin{equation}
\label{eq:local-global}
\HC_*^{I}(M,\alpha)=\bigoplus_x\HC_*(x).
\end{equation}
As a consequence, $\HC_m^{I}(M,\alpha)=0$ for any interval $I$ when
all orbits $x$ with action in $I$ are isolated and have zero local
contact homology in degree $m$. Here, as in Section \ref{sec:ch}, we
can fix a free homotopy class of loops in $M$ or a collection of such
classes.

There are several ways to circumvent the transversality issues in the
construction of the local contact homology. For instance, one can
simply declare \eqref{eq:lch-lfh2} to be the definition of
$\HC_*(x)$. However, the following approach proposed by Michael
Hutchings is from our perspective more aesthetically pleasing. Namely,
when $x=z^k$ where $z$ is simple, there is a natural
$\Z_k$-action by contactomorphisms on the $k$-fold covering of a
neighborhood of $z$. Combined with the continuation maps, this action
gives rise to a $\Z_k$-action on $\HC_*(\tilde{z})$, where $\tilde{z}$
is the $k$-fold covering of $z$, and one can just set
$\HC_*(x)=\HC_*(\tilde{z})^{\Z_k}$. (The composition of a deck
transformation and the continuation map is still $k$-periodic on the
level of homology since continuation maps are canonical.) A similar
construction can be used in the Hamiltonian setting to define
$\HF_*(\varphi)^{\Z_k}$ without assuming the existence of a regular,
one-periodic in time almost complex structure. Then
\eqref{eq:lch-lfh3} becomes a consequence of \eqref{eq:lch-lfh1} and
the definitions.  Note however that with any of these ``indirect''
definitions of the local contact homology, \eqref{eq:local-global} and
other similar facts would become much less obvious.

\subsection{Technicalities} Regardless of the transversality issues,
there are two other problems one has to deal with when defining the
local contact homology, and this is where \eqref{eq:elb} and
\eqref{eq:elb2} can also be useful. (In \cite{HM} a different
approach, at least on the technical level, is used to show that
the local Floer homology is defined and, moreover, well-defined. That
approach is likely to result in somewhat different requirements on the
solutions of the (CR) equation than \eqref{eq:elb} and \eqref{eq:elb2}.)

The first problem is that to have $\p^2=0$ we need to show that the Floer
trajectories for a small non-degenerate perturbation $\alpha'$ of
$\alpha$ asymptotic to some of the orbits $x$ splits into cannot leave
$U$. The $\omega$-energy of such a trajectory $u$ is necessarily small
since $\alpha'$ is close to $\alpha$. Then it follows immediately from
\eqref{eq:elb} that $u$ is confined to~$U$.

The second problem is to show that the resulting homology is
independent of the perturbation $\alpha'$ and whatever auxiliary data
is used. Here again the main point is to localize the trajectories to
a small neighborhood of $x$ by, say, using a variant of
\eqref{eq:elb2}. Once this is done, the proof is identical to the
stability argument for the global or filtered contact homology pointed
out in Section \ref{sec:gen}.

In the context of the holomorphic curves, the proofs of both facts do
not require any new machinery and are in part given above.

\end{document}